\documentclass[]{interact}

\usepackage{epstopdf}
\usepackage{subfigure}

\usepackage[numbers,sort&compress]{natbib}
\bibpunct[, ]{[}{]}{,}{n}{,}{,}

\theoremstyle{plain}
\newtheorem{theorem}{Theorem}[section]
\newtheorem{lemma}[theorem]{Lemma}
\newtheorem{corollary}[theorem]{Corollary}
\newtheorem{proposition}[theorem]{Proposition}

\theoremstyle{definition}

\newtheorem{example}[theorem]{Example}

\theoremstyle{remark}
\newtheorem{remark}{Remark}

\begin{document}

\articletype{ARTICLE TEMPLATE}

\title{Anti-commuting Solutions of the Yang-Baxter-like Matrix Equation}

\author{\name{Mohammed Ahmed Adam Abdalrahman\textsuperscript{a},  Huijian Zhu\textsuperscript{b}, Jiu Ding\textsuperscript{c}} and Qianglian Huang\textsuperscript{a}\thanks{Corresponding
author. Email: huangql@yzu.edu.cn}\\
\affil{ \textsuperscript{a}School of Mathematics, Yangzhou University, Yangzhou 225002, PR China\\
\textsuperscript{b}School of Mathematics and Statistics, Nanfang University, Guangzhou 510970, PR China\\
\textsuperscript{c} School of Mathematics and Natural Sciences, University of Southern Mississippi,  Hattiesburg, MS 39406, USA}
}

\maketitle

\begin{abstract}
\indent We solve the Yang-Baxter-like matrix equation $AXA = XAX$ for a general given matrix $A$ to get
all anti-commuting solutions, by using the Jordan canonical form of $A$ and applying some
new facts on a general homogeneous Sylvester equation. Our main result provides all the
anti-commuting solutions of the nonlinear matrix equation.
\end{abstract}


\section{Introduction}

\quad  In the past decade, motivated by wide applications \cite{ff, yg} of the classic Yang-Baxter
equation, which was independently proposed by Yang \cite{yang} and Baxter \cite{baxter} in 1967 and 1972
respectively, in the mathematical and physical sciences such as knot theory, braid groups, and
quantum groups, solving the {\em Yang-Baxter-like matrix equation}
\begin{eqnarray}
\label{yb}
AXA = XAX,
\end{eqnarray}
which has been so named because of its similarity in format to the original Yang-Baxter equation,
has attracted much attention in the community of linear algebra and other areas such as numerical analysis. In the quadratic equation (\ref{yb}), $A$ is a given $n \times n$ complex matrix and
one desires to find all solutions or some solutions satisfying one or several additional properties.

A solution $B$ of (\ref{yb}) such that $AB = BA$ is called a {\em commuting solution}. All commuting solutions of the Yang-Baxter-like matrix equation for general $A$ have been obtained in
theory via the works of \cite{dd1, dd2, swj} with different approaches. Among them, a particular
class of commuting solutions were explicitly constructed. More specifically, associated with
each eigenvalue of $A$, the product of $A$ and the spectral projection matrix from $\mathbb{C}^n$
onto the generalized eigenspace of the eigenvalue is a commuting solution \cite{dr2}, called the spectral solution with respect to the given eigenvalue. The collection of spectral solutions
has been enlarged to the family of all projection-based commuting solutions in \cite{dd3}.
For the class of diagonalizable matrices $A$, based on a result \cite{ss, xz} on the unique
solution of a Sylvester equation, all commuting solutions have been constructed by \cite{dd}.

As for non-commuting solutions, however, the progress is limited. When $A$ belongs to several special
matrix classes, such as $A$ is an elementary matrix, has a small rank, a small number of
distinct eigenvalues, or a simple minimal polynomial, the solution manifold of (\ref{yb})
has been investigated in the literature; see, for example, \cite{adh, ti, zcd, zd}. But so far,
finding all solutions of (\ref{yb}) for a general matrix $A$ is still an unsolved problem.

Inspired by the techniques employed to find all commuting solutions of the Yang-Baxter-like
matrix equation, in this paper we focus on constructing a class of non-commuting solutions of
(\ref{yb}). Our task is to find all solutions $B$ of (\ref{yb}) satisfying $AB = - BA$.
Such solutions will be referred to as {\em anti-commuting solutions}. The main result,
Theorem 3.9, of the paper will demonstrate that, the problem of obtaining all anti-commuting
solutions for a general matrix $A$ can be reduced to that of finding all anti-commuting
solutions of a simplified Yang-Baxter-like matrix equation of smaller size, which
is determined by the zero eigenvalue of $A$. In particular, if $A$ is nonsingular,
then the zero matrix is the only anti-commuting solution of (\ref{yb}).


Our approach in the following is: To find anti-commuting solutions of the matrix
equation, we solve an equivalent equation to (\ref{yb}) after posing the additional
anti-commutativity condition. The original two matrix equations are simplified by
replacing the given matrix with its Jordan canonical form.

The paper is organized as follows. We prove several preliminary results in the next
section, which are of importance by themselves in addition to their usefulness for
Section 3, in which we first give some equivalent and sufficient formations of the
original problem, followed by considering the special case that $A$ has only one
eigenvalue, and then we apply the results of Section 2 to the general case. Two
numerical examples will be presented in Section 4 to illustrate the main solution
theorem. We conclude with Section 5.

\section{Preliminary Results}

\quad $\;$ In this section we investigate the solution structure of the linear
matrix equation $AX = -XA$ as the first step towards solving for anti-commuting
solutions of the Yang-Baxter-like matrix equation. We shall use Jordan canonical
forms throughout the paper to fulfill our purpose.

The Jordan canonical form $J$ of $A$ is a block diagonal matrix with Jordan blocks as its
diagonal square matrices. Each Jordan block has its main diagonal entries one eigenvalue
of $A$ and the sup-diagonal entries are straight $1$'s, with all other entries $0$.
A $t \times t$ Jordan block with its diagonal entries $\lambda$ is denoted as
$J_t(\lambda)$, which is the following upper triangular matrix:
\[ J_t(\lambda) = \left[ \begin{array}{ccccccccc}
              \lambda & 1       & 0       & \cdots  & 0       & 0      \\
              0       & \lambda & 1       & \cdots  & 0       & 0      \\
              0       & 0       & \lambda & \cdots  & 0       & 0      \\
              \vdots  & \vdots  & \vdots  & \ddots  & \vdots  & \vdots \\
              0       & 0       & 0       & \cdots  & \lambda & 1      \\
              0       & 0       & 0       & \cdots  & 0       & \lambda
              \end{array} \right]. \]
It is well-known that the Jordan block $J_t(0)$ corresponding to eigenvalue $0$
satisfies $J_t(0)^t = 0$, since it is a nilpotent matrix of index $t$.

We need a couple of basic facts about Jordan blocks that are the building bricks of the
Jordan canonical form of a matrix. There may be several Jordan blocks corresponding to
the same eigenvalue $\lambda$ of $A$, the largest size of which is the {\em index} of
$\lambda$ that is, by definition, the smallest positive integer $k$ such that
dim$N[(A - \lambda I)^k] = \mbox{dim}N[(A - \lambda I)^{k+1}]$.

When studying the anti-commutability equality $AB = -BA$, we shall
consider a simplified equality $JK = -KJ$ with $J$ the Jordan canonical form of $A$.
Since $J$ is a block diagonal matrix, we first study a basic equation
\begin{eqnarray}
\label{ba}
J_t(\lambda) Y = - Y J_s(\mu)
\end{eqnarray}
for the $t \times s$ unknown matrix $Y$, with $\lambda$ and $\mu$ complex numbers.
Note that $J_t(\lambda) = \lambda I_t + J_t(0)$ and $J_s(\mu) = \mu I_s + J_s(0)$,
where $I_t$ and $I_s$ are the $t \times t$ and $s \times s$ identity matrix, respectively. Then
\[ J_t(\lambda) Y = [\lambda I_t + J_t(0)]Y = \lambda Y + J_t(0)Y \]
and
\[ Y J_s(\mu) = Y [\mu I_s + J_s(0)] = \mu Y + YJ_s(0). \]
Consequently, a $t \times s$ matrix $K$ satisfies (\ref{ba}) if and only if
\begin{eqnarray}
\label{f}
J_t(0)K = - K[J_s(0) + (\lambda + \mu) I_s].
\end{eqnarray}

The following two lemmas will not only be applied to establishing the main result
for a general Yang-Baxter-like matrix equation in the next section, but they also
provide a basis for a general assertion of this section for solving a homogeneous
Sylvester equation, which extends some classic result (for example, Theorem 5.15
of \cite{cc}) from the commuting case to the anti-commuting one. They deal with the
different cases of $\lambda \neq -\mu$ and $\lambda = - \mu$ for the given eigenvalues
$\lambda$ and $\mu$, respectively.

\begin{lemma}
Let $t$ and $s$ be positive integers and let $\lambda$ and $\mu$ be complex numbers such that
$\lambda \neq -\mu$. Then a $t \times s$ matrix $K$ satisfies the equality
$J_t(\lambda) K = - K J_s(\mu)$ if and only if $K = 0$.
\end{lemma}

\noindent {\bf Proof} Only the necessity part needs a proof. Suppose that $J_t(\lambda) K = -K J_s(\mu)$.
Then (\ref{f}) is satisfied. Since $J_t(0)^t = 0$, using (\ref{f}) repeatedly $t$ times, we see that
\begin{eqnarray*}
0  & = & J_t(0)^t K = J_t(0)^{t-1} \cdot J_t(0) K \\
   & = &  - J_t(0)^{t-1} K [J_s(0) + (\lambda + \mu)I_s] \\
   & = & - J_t(0)^{t-2} \cdot J_t(0) K \cdot [J_s(0) + (\lambda + \mu)I_s] \\
   & = & J_t(0)^{t-2} K [J_s(0) + (\lambda + \mu)I_s]^2 \\
   & = & \cdots = (-1)^{t-1} \cdot J_t(0) K \cdot [J_s(0) + (\lambda + \mu)I_s]^{t-1} \\
   & = & (-1)^t K [J_s(0) + (\lambda + \mu)I_s]^t.
\end{eqnarray*}
The assumption $\lambda \neq -\mu$ guarantees that $J_s(0) + (\lambda + \mu)I_s$ is a nonsingular matrix, so the
equality $0 = (-1)^t K [J_s(0) + (\lambda + \mu)I_s]^t$ implies $K = 0$. \hfill $\Box$

\begin{corollary}
Let $\lambda \neq 0$. A $t \times s$ matrix $K$ solves the equation $J_t(\lambda) Y = - Y J_s(\lambda)$
if and only if $K = 0$.
\end{corollary}

\begin{lemma}
Let $t$ and $s$ be positive integers and let $\lambda$ and $\mu$ be complex numbers such that
$\lambda = -\mu$. Then a $t \times s$ matrix $K$ satisfies the equality
$J_t(\lambda) K = - K J_s(\mu)$ if and only if $K = [0 \; \; \hat{K}]$ when $t \le s$ or
\[ K = \left[ \begin{array}{c}
                \hat{K} \\
                  0
               \end{array} \right] \]
when $t \ge s$, where $\hat{K}$ is an upper triangular matrix of the form
\begin{eqnarray}
\label{hatK}
\hat{K} = \left[ \begin{array}{rrrrrr}
\hat{k}_1 & \hat{k}_2  & \hat{k}_3  & \cdots & \hat{k}_{r-1}       & \hat{k}_r           \\
0         & -\hat{k}_1 & -\hat{k}_2 & \cdots & -\hat{k}_{r-2}      & -\hat{k}_{r-1}      \\
0         & 0          & \hat{k}_1  & \cdots & \hat{k}_{r-3}       & \hat{k}_{r-2}       \\
\vdots    & \vdots     & \vdots     & \ddots & \vdots              & \vdots              \\
0         & 0          & 0          & \cdots & (-1)^{r-2}\hat{k}_1 & (-1)^{r-2}\hat{k}_2 \\
0         & 0          & 0          & \cdots & 0                   & (-1)^{r-1}\hat{k}_1
                \end{array} \right],
\end{eqnarray}
with $r = \min\{t, s\}$ and $\hat{k}_1, \ldots, \hat{k}_r$ are arbitrary complex numbers.
\end{lemma}

\noindent {\bf Proof} We only prove the lemma for the case of $t \le s$; the other case can be done
by the same idea. Since $\lambda + \mu = 0$, it is clear that $J_t(\lambda) K = -K J_s(\mu)$
if and only if
\begin{eqnarray}
\label{zero}
J_t(0) K = -K J_s(0).
\end{eqnarray}

First we verify that, if $K = [0 \; \; \hat{K}]$ with $\hat{K}$ given by (\ref{hatK}), then
(\ref{zero}) is true. This follows from the observation that
\[ J_t(0)K = J_t(0) [0 \; \; \hat{K}] = [0 \; \; J_t(0) \hat{K}] \]
and
\[ -KJ_s(0) = -[0 \; \; \hat{K}] J_s(0) = [0 \; \; -\hat{K} J_t(0)], \]
and the direct computation:
\begin{eqnarray*}
&   & J_t(0) \hat{K} \\
& = & \left[ \begin{array}{ccccccccc}
              0       & 1      & 0      & \cdots  & 0      & 0      \\
              0       & 0      & 1      & \cdots  & 0      & 0      \\
              0       & 0      & 0      & \cdots  & 0      & 0      \\
              \vdots  & \vdots & \vdots & \ddots  & \vdots & \vdots \\
              0       & 0      & 0      & \cdots  & 0      & 1      \\
              0       & 0      & 0      & \cdots  & 0      & 0
              \end{array} \right] \left[ \begin{array}{rrrrrr}
 \hat{k}_1 & \hat{k}_2  & \hat{k}_3  & \cdots & \hat{k}_{r-1}       & \hat{k}_r           \\
 0         & -\hat{k}_1 & -\hat{k}_2 & \cdots & -\hat{k}_{r-2}      & -\hat{k}_{r-1}      \\
 0         & 0          & \hat{k}_1  & \cdots & \hat{k}_{r-3}       & \hat{k}_{r-2}       \\
 \vdots    & \vdots     & \vdots     & \ddots & \vdots              & \vdots              \\
 0         & 0          & 0          & \cdots & (-1)^{r-2}\hat{k}_1 & (-1)^{r-2}\hat{k}_2 \\
 0         & 0          & 0          & \cdots & 0                   & (-1)^{r-1}\hat{k}_1
                                         \end{array} \right]                              \\
& = & \left[ \begin{array}{rrrrrr}
              0      & -\hat{k}_1 & -\hat{k}_2 & \cdots & -\hat{k}_{t-2} & -\hat{k}_{t-1}      \\
              0      & 0          & \hat{k}_1  & \cdots & \hat{k}_{t-3}  & \hat{k}_{t-2}       \\
              0      & 0          & 0          & \cdots & -\hat{k}_{t-4} & -\hat{k}_{t-3}      \\
              \vdots & \vdots     & \vdots     & \ddots & \vdots         & \vdots              \\
              0      & 0          & 0          & \ddots & 0              & (-1)^{t-1}\hat{k}_1 \\
              0      & 0          & 0          & \cdots & 0              & 0
             \end{array} \right],
\end{eqnarray*}
{\small
\begin{eqnarray*}
&   & -\hat{K} J_t(0) \\
& = & -\left[ \begin{array}{rrrrrr}
\hat{k}_1 & \hat{k}_2  & \hat{k}_3  & \cdots & \hat{k}_{r-1}       & \hat{k}_r           \\
0         & -\hat{k}_1 & -\hat{k}_2 & \cdots & -\hat{k}_{r-2}      & -\hat{k}_{r-1}      \\
0         & 0          & \hat{k}_1  & \cdots & \hat{k}_{r-3}       & \hat{k}_{r-2}       \\
\vdots    & \vdots     & \vdots     & \ddots & \vdots              & \vdots              \\
0         & 0          & 0          & \cdots & (-1)^{r-2}\hat{k}_1 & (-1)^{r-2}\hat{k}_2 \\
0         & 0          & 0          & \cdots & 0                   & (-1)^{r-1}\hat{k}_1
              \end{array} \right] \left[ \begin{array}{ccccccccc}
              0       & 1      & 0      & \cdots  & 0      & 0      \\
              0       & 0      & 1      & \cdots  & 0      & 0      \\
              0       & 0      & 0      & \cdots  & 0      & 0      \\
              \vdots  & \vdots & \vdots & \ddots  & \vdots & \vdots \\
              0       & 0      & 0      & \cdots  & 0      & 1      \\
              0       & 0      & 0      & \cdots  & 0      & 0
              \end{array} \right] \\
& = & \left[ \begin{array}{rrrrrr}
              0      & -\hat{k}_1 & -\hat{k}_2 & \cdots & -\hat{k}_{t-2} & -\hat{k}_{t-1}      \\
              0      & 0          & \hat{k}_1  & \cdots & \hat{k}_{t-3}  & \hat{k}_{t-2}       \\
              0      & 0          & 0          & \cdots & -\hat{k}_{t-4} & -\hat{k}_{t-3}      \\
              \vdots & \vdots     & \vdots     & \ddots & \vdots         & \vdots              \\
              0      & 0          & 0          & \ddots & 0              & (-1)^{t-1}\hat{k}_1 \\
              0      & 0          & 0          & \cdots & 0              & 0
             \end{array} \right] = J_t(0) \hat{K}.
\end{eqnarray*}}

Conversely, let a $t \times s$ matrix $K = [k_{ij}]$ satisfy (\ref{zero}).
Because of the special structure of $J_t(0)$, we see that rows $1$
through $t-1$ of $J_t(0) K$ are rows $2$ through $t$ of $K$ respectively and row $t$
of $J_t(0) K$ is zero. Similarly, column $1$ of $K J_s(0)$ is zero and columns $2$
through $s$ of $KJ_s(0)$ are columns $1$ through $s-1$ of $K$ respectively. Thus, $J_t(0) K = -K J_s(0)$
implies
\[ \left[ \begin{array}{lllll}
              k_{21} & k_{22} & \cdots & k_{2,s-1} & k_{2s} \\
              k_{31} & k_{32} & \cdots & k_{3,s-1} & k_{3s} \\
              \vdots & \vdots & \ddots & \vdots    & \vdots \\
              k_{t1} & k_{t2} & \cdots & k_{t,s-1} & k_{ts} \\
              0      & 0      & \cdots & 0             & 0
             \end{array} \right] = - \left[ \begin{array}{lllll}
              0      & k_{11}    & \cdots & k_{1,s-2}   & k_{1,s-1}   \\
              0      & k_{21}    & \cdots & k_{2,s-2}   & k_{2,s-1}   \\
              \vdots & \vdots    & \ddots & \vdots      & \vdots      \\
              0      & k_{t-1,1} & \ddots & k_{t-1,s-2} & k_{t-1,s-1} \\
              0      & k_{t1}    & \cdots & k_{t,s-2}   & k_{t,s-1}
             \end{array} \right]. \]
It follows that $k_{i1} = 0$ for $i = 2, \ldots, t$ and $k_{tj} = 0$ for $j = 1, \ldots, s-1$.
Comparing the corresponding entries of the both side matrices, we see that
\[ k_{i+1, j} = - k_{i, j-1}, \; \; i = 1, \ldots, t-1, \; j = 2, \ldots, s, \]
The above recursive relation means that when the row and column indices increase by $1$ each,
the corresponding entry changes the sign. Hence, since all the entries of row $t$ of $K$ are $0$
except for the last one, starting from $k_{t-1,2} = - k_{t3} = 0$, we obtain in succession that
$k_{ij} = 0$ for all $(i, j)$ with $1 \le i \le t$ and $1 \le j \le s - t$, so $K = [0 \; \; \hat{K}]$
with $0$ the $t \times (s - t)$ zero matrix and $\hat{K} = [\hat{k}_{ij}]$ is $t \times t$. The same
argument shows that $\hat{k}_{ij} = 0$ for all $i, j = 1, \ldots, t$ such that $i > j$ and
$\hat{k}_{ii} = (-1)^{i-1}\hat{k}_{11}$ for $i = 2, \ldots, t$, $\hat{k}_{i,i+1} =
(-1)^{i-1}\hat{k}_{12}$ for $i = 2, \ldots, t-1$, \ldots, $\hat{k}_{2t} = -\hat{k}_{1, t-1}$.
Hence $\hat{K}$ is given by (\ref{hatK}) with $\hat{k}_i \equiv \hat{k}_{1i}$ being arbitrary
numbers for $i = 1, \ldots t$.  \hfill $\Box$

\begin{corollary}
In particular, Lemma 2.3 gives all solutions of the equation $J_t(0) Y = - Y J_s(0)$.
\end{corollary}

Denote by $J(\lambda)$ a block diagonal matrix with its diagonal consisting of several
Jordan blocks associated to the same eigenvalue $\lambda$. Then Lemmas 2.1 and 2.3 can
be extended to the following proposition.

\begin{proposition}
Let $K$ be a solution of the matrix equation
\begin{eqnarray*}
J(\lambda) Y = - Y J(\mu),
\end{eqnarray*}
where $J(\lambda) = \mbox{\em diag} [J_{r_1}(\lambda), \ldots, J_{r_m}(\lambda)]$ and $J(\mu) =
\mbox{\em diag} [J_{s_1}(\mu), \ldots, J_{s_l}(\mu)]$.

{\em (i)} If $\lambda \neq -\mu$, then $K = 0$.

{\em (ii)} If $\lambda = - \mu$, then $K$ is an $m \times l$ block matrix with its
$(i, j)$ block an $r_i \times s_j$ matrix $K_{ij}$ that has the same pattern as given in Lemma 2.3.
\end{proposition}

\noindent {\bf Proof} Partition $K$ according to the diagonal blocks size of $J(\lambda)$ and
$J(\mu)$, which is consistent with the equality $J(\lambda) K = - K J(\mu)$. Using the block
matrix multiplication, we have
\[ J_{r_i}(\lambda) K_{ij} = - K_{ij} J_{s_j}(\mu), \; \; i = 1, \ldots, m, \; j = 1, \ldots, l. \]
Then immediately, Lemma 2.1 implies (i) and Lemma 2.3 gives (ii). \hfill $\Box$

We generalize the above proposition to the case of multi-eigenvalues. For this purpose,
consider a general linear matrix equation
\begin{eqnarray}
\label{lme}
U X = - X V,
\end{eqnarray}
where $U$ and $V$ are respectively $u \times u$ and $v \times v$ known matrices,
so the unknown matrix $X$ is $u \times v$.

Let $J_U$ and $J_V$ be the Jordan canonical forms of $U$ and $V$, respectively.
There are nonsingular matrices $P$ and $Q$ such that $U = P J_U P^{-1}$ and $V = Q J_V Q^{-1}$.
Then the equation (\ref{lme}) for $X$ is equivalent to the equation
\begin{eqnarray}
\label{lme2}
J_U Y = - Y J_V
\end{eqnarray}
for $Y$ with $Y = P^{-1} X Q$. Write
\[ J_U = \mbox{diag} [J(\lambda_1), \ldots, J(\lambda_p)] \; \; \; \mbox{and} \; \; \;
 J_V = \mbox{diag} [J(\mu_1), \ldots, J(\mu_q)], \]
where $\lambda_1, \ldots, \lambda_p$ are all distinct eigenvalues of $U$ and
$\mu_1, \ldots, \mu_q$ are all distinct eigenvalues of $V$.
Partitioning $Y = [Y_{ij}]$ accordingly as a $p \times q$ block matrix,
we see that solving (\ref{lme2}) is equivalent to solving the $pq$ sub-matrix equations
\begin{eqnarray*}
J(\lambda_i) Y_{ij} = - Y_{ij} J(\mu_j), \; \; i = 1, \ldots, p; \; j = 1, \ldots, q.
\end{eqnarray*}

Proposition 2.5 then gives rise to the main result of this section.

\begin{theorem}
The solutions of $U X = - X V$ are $X = P Y Q^{-1}$ in which the block matrix $Y = [Y_{ij}]$
has the property: For $i = 1, \ldots, p$ and $j = 1, \ldots, q$,

{\em (i)} if $\lambda_i \neq -\mu_j$, then $Y_{ij} = 0$;

{\em (ii)} if $\lambda_i = -\mu_j$, then $Y_{ij}$ can be partitioned as a block matrix with
each block as given by {\em (ii)} in Proposition 2.5.
\end{theorem}



\begin{corollary}
If the eigenvalues $\lambda_1, \ldots, \lambda_p$ of $U$ and $\mu_1, \ldots, \mu_q$ of $V$
satisfy the condition that $\lambda_i \neq -\mu_j$ for all $i = 1, \ldots, p$ and $j = 1,
\ldots, q$, then $X = 0$ is the only solution of the equation $U X = - X V$.
\end{corollary}

\begin{corollary}
If the eigenvalues $\lambda_1, \ldots, \lambda_p$ of $U$ are nonzero and
$\lambda_i \neq -\lambda_j$ for all $i \neq j$, then $X = 0$ is the only solution
of the equation $U X = - X U$.
\end{corollary}

\section{Anti-commuting Solutions of the Yang-Baxter-like Matrix Equation}

\quad Suppose that $A$ is a given $n \times n$ complex matrix. We first give some useful
results that will be applied to proving the main theorems of the paper. The first lemma
below provides an equivalent condition for a matrix to be an anti-commuting solution of the Yang-Baxter-like matrix equation.

\begin{lemma}
Let an $n \times n$ complex matrix $B$ satisfy the equality $AB = -BA$. Then $B$ solves (\ref{yb}) if and only if
\begin{eqnarray}
\label{ec1}
 B(B - A)A = 0.
\end{eqnarray}
\end{lemma}

\noindent {\bf Proof} Since $AB = -BA$,
\[ ABA - BAB = -BAA + BBA = B(B - A)A. \]
Thus, $ABA = BAB$ if and only if $B(B - A)A = 0$. \hfill $\Box$

\begin{remark} Another equivalent condition for $B$ to be an anti-commuting solution of
(\ref{yb}) is $A(B - A)B = 0$, which follows from $ABA - BAB = -AAB + ABB = A(B - A)B$.
\end{remark}

\begin{remark}
Suppose that $AB = -BA$. Then a sufficient condition for (\ref{ec1}) to be valid is $BA = A^2$
or $BA = B^2$, since the former implies $(B - A)A = 0$ and the latter implies $B(B - A) = 0$,
from which $B(B - A)A = 0$.
\end{remark}



It is well-known (see, for example, Lemma 3.1 of \cite{drz1}) that, if $J$ is the Jordan
canonical form of $A$, then solving (\ref{yb}) for $X$ is equivalent to solving the
simplified Yang-Baxter-like matrix equation
\begin{eqnarray}
\label{yb1}
JYJ = YJY
\end{eqnarray}
for $Y$ in the sense that $B$ is a solution of (\ref{yb}) if and only if $K$ is a solution
of (\ref{yb1}) with their relation  $B = W K W^{-1}$, where $W$ is the similarity matrix
connecting the given matrix $A$ to its Jordan canonical form $J$, that is, $A = W J W^{-1}$. Furthermore, $B$ and $K$ share some common properties, such as $B$ commutes with $A$ if and
only if $K$ commutes with $J$.

The same idea implies that with the similarity relation $B = W K W^{-1}$,
the matrix $B$ is an anti-commuting solution of (\ref{yb}) if and only if the
matrix $K$ is an anti-commuting solution of (\ref{yb1}). Thus in our analysis
below, it is sufficient to find all the anti-commuting solutions of (\ref{yb1}).
As a consequence, all anti-commuting solutions in this section will be stated
with respect to the simplified Yang-Baxter-like matrix equation (\ref{yb1}),
from which all anti-commuting solutions of the original Yang-Baxter-like matrix
equation (\ref{yb}) are available immediately via the similarity matrix $W$.


Our strategy is, first we find all $n \times n$ matrices $K$ satisfying the
anti-commutability equation $JY = -YJ$, for which we can apply the general results
of the previous section; then, among all such matrices $K$ we search for those
satisfying the equation $JYJ = YJY$. By Lemma 3.1, this latter task is reduced
to solving the homogeneous equation $Y(Y-J)J = 0$.

We first consider the simplest case that $J = J_n(\lambda)$, which can be analyzed
by applying Corollaries 2.2 and 2.4.

\begin{theorem}
If $\lambda\neq 0$, then $K=0$ is the only anti-commuting solution of
$J_n(\lambda)YJ_n(\lambda)=YJ_n(\lambda)Y$. If $\lambda=0$, then all
anti-commuting solutions of $J_n(0)YJ_n(0)=YJ_n(0)Y$ are
\[ K = \left[ \begin{array}{cc}
              0 & x    \\
              0 & 0
\end{array} \right] \; \; \mbox{for} \; \; n = 2, \]
\[ K = \left[ \begin{array}{rrr}
              0 & y & x  \\
              0 & 0 & -y \\
              0 & 0 & 0
\end{array} \right] \; \; \mbox{for} \; \; n = 3, \; \; \mbox{and} \]
\begin{equation}
\label{K31}
K = \left[ \begin{array}{rrrrrr}
              0      & 0      & \cdots & 0      & y      & x      \\
              0      & 0      & \cdots & 0      & 0      & -y     \\
              0      & 0      & \cdots & 0      & 0      & 0      \\
              \vdots & \vdots & \ddots & \vdots & \vdots & \vdots \\
              0      & 0      & 0      & \cdots & 0      & 0      \\
              0      & 0      & 0      & \cdots & 0      & 0
\end{array} \right],
\end{equation}
where $x$ and $y$ are arbitrary complex numbers.
\end{theorem}

\noindent {\bf Proof} The result for $\lambda \neq 0$ is immediate from Corollary 2.2
applied to $J_n(\lambda)Y = -YJ_n(\lambda)$. Suppose that $\lambda = 0$. The case
$n \le 3$ is easy to verify, so we assume that $n \ge 4$. By Corollary 2.4 with
$t = s = n$, the condition $J_n(0) K = - K J_n(0)$ implies that $K$ has the
structure given by (\ref{hatK}). In other words, all such $K$ are
\begin{eqnarray}
\label{K}
K = \left[ \begin{array}{rrrrrr}
              k_1    & k_2    & k_3    & \cdots & k_{n-1}       & k_n           \\
              0      & -k_1   & -k_2   & \cdots & -k_{n-2}      & -k_{n-1}      \\
              0      & 0      & k_1    & \cdots & k_{n-3}       & k_{n-2}       \\
              \vdots & \vdots & \vdots & \ddots & \vdots        & \vdots        \\
              0      & 0      & 0      & \cdots & (-1)^{n-2}k_1 & (-1)^{n-2}k_2 \\
              0      & 0      & 0      & \cdots & 0             & (-1)^{n-1}k_1
             \end{array} \right],
\end{eqnarray}
in which $k_1, \ldots, k_n$ are arbitrary complex numbers.

The remaining thing is to select $K$ among the above matrices (\ref{K}) to
satisfy the second requirement $K[K - J_n(0)]J_n(0) = 0$. A direct computation gives that
$[K - J_n(0)]J_n(0) =$
\[ \left[ \begin{array}{rrrrrrr}
 0      & k_1    & k_2 - 1 & k_3      & \cdots & k_{n-2}       & k_{n-1}           \\
 0      & 0      & -k_1    & -k_2 - 1 & \cdots & -k_{n-3}      & -k_{n-2}          \\
 0      & 0      & 0       & k_1      & \cdots & k_{n-4}       & k_{n-3}           \\
 \vdots & \vdots & \vdots  & \vdots   & \ddots & \vdots        & \vdots            \\
 0      & 0      & 0       & 0        & \cdots & (-1)^{n-3}k_1 & (-1)^{n-3}k_2 - 1 \\
 0      & 0      & 0       & 0        & \cdots & 0             & (-1)^{n-2}k_1     \\
 0      & 0      & 0       & 0        & \cdots & 0             & 0
             \end{array} \right], \]
so the entries of the first row of $K[K - J_n(0)]J_n(0)$ is
$0, k_1^2, - k_1, 2k_1k_3 - k_2(k_2 + 1), \ldots, k_1k_{n-1} - k_2k_{n-2}
+ \cdots + k_{n-2}[(-1)^{n-3}k_2 - 1] + (-1)^{n-2}k_{n-1}k_1$.
Since the $(1, 2)-$entry and $(1, 3)-$entry of $K[K - J_n(0)]J_n(0) = 0$
are zero, $k_1 = 0$ and $k_2 = 0$ or $-1$. But since the $(2, 5)-$entry of
$K[K - J_n(0)]J_n(0)$ is zero, the possibility of $k_2 = -1$ is excluded,
so $k_2 = 0$. Then we obtain in succession that $k_3 = k_4 = \cdots =
k_{n-2} = 0$, with $k_{n-1}$ and $k_n$ being arbitrary numbers. Hence all
the solutions $K$ of (\ref{yb1}) are given by (\ref{K31}). \hfill $\Box$

Now we let $A$ be a general $n \times n$ complex matrix. For the simplicity
of presentation and analysis, we can write the Jordan form $J$ of $A$ as
\begin{eqnarray}
\label{J}
J = \left[ \begin{array}{cccc}
       J(\lambda_1) &              &        &     \\
                    & J(\lambda_2) &        &     \\
                    &              & \ddots &     \\
                    &     &        & J(\lambda_d)
              \end{array} \right],
\end{eqnarray}
in which $\lambda_1, \ldots, \lambda_d$ are all distinct eigenvalues of $A$ and
\[ J(\lambda_i) = \left[ \begin{array}{ccc}
       J_{t_{i1}}(\lambda_i) &        &        \\
                          & \ddots &        \\
                          &        & J_{t_{im_i}}(\lambda_i)
              \end{array} \right], \; \; i = 1, \ldots, d. \]

In the following, we first consider the case that $J$ is a nonsingular matrix,
so its eigenvalues are all nonzero.

\begin{theorem} If $J$ is nonsingular, then $K = 0$ is the only anti-commuting
solution of $JYJ = YJY$.
\end{theorem}

\noindent {\bf Proof} Let the $d$ distinct eigenvalues of $J$, which are all
nonzero by assumption, be written as
\[ \lambda_1, \ldots, \lambda_{d-2k}, \mu_1, -\mu_1, \ldots, \mu_k, -\mu_k, \]
in which $-\lambda_i$ is not an eigenvalue of $J$ for $i = 1, \ldots, d-2k$ with
$0 \le k \le d/2$. Without loss of generality, we assume that
\begin{eqnarray}
\label{JJ}
J = \left[ \begin{array}{cccc}
       J_1 &   0   \\
         0 & J_2
           \end{array} \right],
           \end{eqnarray}
where $J_1$ and $J_2$ correspond to the eigenvalues $\lambda_1, \ldots, \lambda_{d-2k},
\mu_1, \ldots, \mu_k$ and $-\mu_1, \ldots, -\mu_k$, respectively.

Suppose that $K$ is an anti-commuting solution of (\ref{yb1}).
Partition $K$ as a $2 \times 2$ block matrix
\begin{eqnarray*}
\label{K2}
K = \left[ \begin{array}{ccc}
       K_1 & K_2 \\
       K_3 & K_4
           \end{array} \right]
\end{eqnarray*}
according to the sizes of the diagonal blocks of $J$ as in (\ref{JJ}).
Then the equality $JK = -KJ$ becomes
\[ \left[ \begin{array}{lllll}
              J_1K_1 & J_1K_2 \\
              J_2K_3 & J_2K_4
             \end{array} \right] = - \left[ \begin{array}{lllll}
             K_1J_1 & K_2J_2  \\
             K_3J_1 & K_4J_2
             \end{array} \right]. \]
By Corollary 2.8, $K_1 = 0$ from $J_1K_1 = -K_1 J_1$ and $K_4 = 0$ from
$J_2K_4 = -K_4 J_2$. On the other hand, Lemma 3.1 ensures $K(K - J)J = 0$,
from which $KJ = K^2$ since $J$ is nonsingular. It follows that
\[ \left[ \begin{array}{cc}
               0        & K_2J_2 \\
               K_3J_1   & 0
             \end{array} \right] = \left[ \begin{array}{cc}
             K_2K_3 &  0   \\
             0      & K_3K_2
             \end{array} \right]. \]
Since $J_1$ and $J_2$ are nonsingular, $K_2J_2 = 0$ implies $K_2 = 0$ and
$K_3J_1 = 0$ implies $K_3 = 0$. Therefore, $K= 0$. \hfill $\Box$

Now we consider the remaining case that the given matrix $A$ is singular, in
other words, $0$ is an eigenvalue of $A$. Due to the fact that $J$ is a block diagonal matrix, we need the following lemma that
concerns the structure of $K$ in the equality $HK = - KH$ when $H$ is a block diagonal
matrix.

\begin{lemma}
Let $H = \mbox{\em diag} (H_1, \ldots, H_d)$ with $H_i$ being $h_i \times h_i$ for
$i = 1, \ldots, d$. Suppose that $HK = -KH$ for a $d \times d$ block matrix
$K = [K_{ij}]$ with $K_{ij}$ being $h_i \times h_j$ for $i, j = 1, \ldots, d$.
Then for each pair $(i, j)$, if no eigenvalue of $H_i$ is opposite to any
eigenvalue of $H_j$ in sign, then $K_{ij} = 0$; otherwise, $K_{ij} \neq 0$.
\end{lemma}

\noindent {\bf Proof} Multiplying
\[ \mbox{diag} (H_1, \ldots, H_d) [K_{ij}] = - [K_{ij}] \mbox{diag} (H_1, \ldots, H_d) \]
out gives
\[ H_i K_{ij} = - K_{ij} H_j, \; \; \forall \; i, j = 1, \ldots, d. \]
For the given $i$ and $j$, since $H_i$ and $H_j$ have no eigenvalues
with opposite signs, Theorem 2.6 (i) guarantees that $K_{ij} = 0$. Otherwise,
there are at least one eigenvalue $\lambda$ of $H_i$ and one eigenvalue $\mu$ of
$H_j$ such that $\lambda = -\mu$, Theorem 2.6 (ii) implies that, as a block
matrix, all blocks of $K_{ij}$ are nonzero matrices, which are structured as
in (ii) of Proposition 2.5. \hfill $\Box$

\begin{corollary}
Under the conditions of Lemma 3.2, if in addition no eigenvalue of $H_i$ is
opposite in sign to any eigenvalue of $H_j$ for all $i \neq j$, then $K$ is a block
diagonal matrix.
\end{corollary}

\noindent {\bf Proof} By Corollary 2.7, $K = \mbox{diag} (K_1, \ldots, K_d)$ with
$K_i = K_{ii}$ for $i = 1, \ldots, d$. \hfill $\Box$


Now we examine the equation $JY = -YJ$. According to the block structure of $J$
as in (\ref{J}), partition $K$ as
\begin{eqnarray}
\label{KB}
K = \left[ \begin{array}{ccc}
       K_{11} & \cdots & K_{1d} \\
       \vdots & \vdots & \vdots \\
       K_{d1} & \cdots & K_{dd}
              \end{array} \right].
\end{eqnarray}
Then $JK= -KJ$ if and only if
\begin{eqnarray}
\label{kij}
J(\lambda_i) K_{ij} = - K_{ij} J(\lambda_j), \; i, j = 1, \ldots, d.
\end{eqnarray}

Proposition 2.5 directly gives rise to the following result.

\begin{proposition}
Let $K = [K_{ij}]$ be a solution of the equation $JY = -YJ$.
Then for each pair $(i, j)$, if $\lambda_i \neq -\lambda_j$, then $K_{ij} = 0$;
if $\lambda_i = -\lambda_j$, then $K_{ij}$ is an $m_i \times m_j$ block matrix
with its $(k, l)$ block a $t_{ik} \times t_{il}$ matrix whose pattern is
given in Lemma 2.3.
\end{proposition}

Based on the above results, the next theorem gives the structure of all
anti-commuting solutions of (\ref{yb1}) for a general Jordan canonical form
$J$ of $A$.

\begin{theorem}
Let the distinct eigenvalues $\lambda_1, \ldots, \lambda_d$ of $J$ be such
that $\lambda_i \neq - \lambda_j$ for all $i \neq j$.
Then all anti-commuting solutions of $JYJ = YJY$ are
$K = \mbox{\em diag}(K_1, \ldots, K_d)$ in which for $i = 1, \ldots, d$,

{\em (i)} if $\lambda_i \neq 0$, then $K_i = 0$;

{\em (ii)} if $\lambda_i = 0$, then $K_i$ is an $m_i \times m_i$ block matrix
with each block structured in Lemma 2.3 as stated by Corollary 2.4 and
satisfies $K_i[K_i - J(0)]J(0) = 0$.
\end{theorem}

\noindent {\bf Proof} The given assumption on the eigenvalues of $A$, by Corollary 3.5,
ensures that $K_{ij} = 0$ is the only solution of (\ref{kij}) whenever $i \neq j$. Hence
all anti-commuting solutions $K$ of $JY = -YJ$ are $d \times d$ block diagonal matrices.
Conclusions (i) and (ii) are immediate from Proposition 3.1. \hfill $\Box$

In Theorem 3.7 (ii), the structure of $K_i$ corresponding to $\lambda_i = 0$
is determined by the equality $K_i[K_i - J(0)]J(0) = 0$, which can be
decomposed into a system of $m_i^2$ matrix equations according to the diagonal
block structure of $J(0)$. As an application of Theorem 3.7, we consider a
special case that $J$ is singular and the eigenvalue $0$ possesses only one
Jordan block. The truth of the following result is then obvious by means of
Theorems 3.2 and 3.7.

\begin{corollary}
Let $n \ge 4$. Suppose that all the distinct eigenvalues of $J$ are
$\lambda_1 = 0, \lambda_2, \ldots, \lambda_d$ with $J(0) = J_m(0)$ and
$\lambda_i \neq -\lambda_j$ for $i \neq j$ from $2$ to $d$. Then all anti-commuting
solutions of the equation $JYJ = YJY$ are $K = \mbox{\em diag}(K_1, 0, \ldots, 0)$,
in which $K_1$ is an $m \times m$ matrix whose expression is given by (\ref{K31}).
\end{corollary}

Using the same method for the proof of Theorem 3.2, we can get rid of the additional
assumption that $\lambda_i \neq -\lambda_j$ for $i \neq j$ from $2$ to $d$ in Theorem 3.7
and Corollary 3.8.
For this purpose, we write $J(0)$ as $J_0$ and rewrite all the nonzero
eigenvalues of $J$ as
\[ \mu_1, \ldots, \mu_{d-2k-1}, \nu_1, -\nu_1, \ldots, \nu_k, -\nu_k, \]
in which $-\mu_i$ is not an eigenvalue of $J$ for $i = 1, \ldots, d-2k-1$ with
$0 \le k < d/2$. Then
\begin{eqnarray}
\label{JJJJ}
J = \left[ \begin{array}{cccc}
  J_0 & 0   &   0   \\
    0 & J_1 &   0   \\
    0 & 0   & J_2
              \end{array} \right],
\end{eqnarray}
where $J_1$ and $J_2$ correspond to the eigenvalues $\mu_1, \ldots, \mu_{d-2k-1},
\nu_1, \ldots, \nu_k$ and $-\nu_1, \ldots, -\nu_k$, respectively.

In $JK = -KJ$, partitioning $K$ as a $3 \times 3$ block matrix
\[ K = \left[ \begin{array}{ccc}
		K_1 & K_2 & K_3 \\
		K_4 & K_5 & K_6 \\
		K_7 & K_8 & K_9
              \end{array} \right] \]
according to the block sizes of $J$ by (\ref{JJJJ}), we obtain that
\[ \left[ \begin{array}{lllll}
             J_0 K_1 & J_0 K_2 & J_0 K_3 \\
             J_1 K_4 & J_1 K_5 & J_1 K_6 \\
             J_2 K_7 & J_2 K_8 & J_2 K_9
             \end{array} \right] = - \left[ \begin{array}{lllll}
             K_1 J_0 & K_2 J_1 & K_3 J_2   \\
             K_4 J_0 & K_5 J_1 & K_6 J_2   \\
             K_7 J_0 & K_8 J_1 & K_9 J_2
             \end{array} \right]. \]
From Corollary 2.3,  $J_0 K_2 = -K_2 J_1$ implies $K_2 = 0, J_0 K_3 = -K_3 J_2$
implies $K_3 = 0, J_1 K_4 = -K_4 J_0$ implies $K_4 = 0, J_1 K_5 = -K_5 J_1$ implies
$K_5 = 0, J_2 K_7 = -K_7 J_0$ implies $K_7 = 0$, and $J_2 K_9 = -K_9 J_2$ implies
$K_9 = 0$. Thus
\[ K = \left[ \begin{array}{ccc}
		K_1 & 0 &   0 \\
		  0 & 0 & K_6 \\
		  0 & K_8 & 0
              \end{array} \right]. \]
Lemma 3.1 gives rise to $K(K - J)J = 0$, which is
\begin{eqnarray*}
 &  & \left[ \begin{array}{ccc}
		K_1 & 0 &   0 \\
		  0 & 0 & K_6 \\
		  0 & K_8 & 0
              \end{array} \right] \left[ \begin{array}{ccc}
		K_1 - J_0 & 0    &   0  \\
		  0       & -J_1 & K_6  \\
		  0       & K_8  & -J_2
              \end{array} \right] \left[ \begin{array}{cccc}
  J_0 & 0   &   0   \\
    0 & J_1 &   0   \\
    0 & 0   & J_2
              \end{array} \right] \\
  & = & \left[ \begin{array}{ccc}
		K_1 & 0 &   0 \\
		  0 & 0 & K_6 \\
		  0 & K_8 & 0
              \end{array} \right] \left[ \begin{array}{ccc}
		(K_1 - J_0)J_0 & 0  &   0     \\
		  0       & -J_1^2  & K_6 J_2 \\
		  0       & K_8J_1  & -J_2^2
              \end{array} \right] \\
  & = & \left[ \begin{array}{ccc}
		K_1(K_1 - J_0)J_0 & 0             &   0        \\
		  0               & K_6K_8J_1     & -K_6 J_2^2 \\
		  0               & -K_8J_1^2     & K_8 K_6J_2
              \end{array} \right] = \left[ \begin{array}{ccc}
		0 & 0 & 0 \\
		0 & 0 & 0 \\
		0 & 0 & 0
              \end{array} \right].
\end{eqnarray*}
Hence $K_1(K_1 - J_0)J_0 = 0, K_6K_8J_1 = 0, -K_6 J_2^2 = 0,
-K_8J_1^2 = 0$, and $K_8K_6J_1 = 0$. Since $J_1$ and $J_2$ are nonsingular,
$K_6 = 0$ from $-K_6J_2^2 = 0$ and $K_8 = 0$ from $-K_8J_1^2 = 0$.
Therefore we have the following theorem.


\begin{theorem}
Let the Jordan canonical form $J$ of $A$ be given by (\ref{JJJJ}). Then all
anti-commuting solutions of $JYJ = YJY$ are
\[ K = \left[ \begin{array}{cccc}
  K_1 &  0  &  0  \\
  0   &  0  &  0  \\
  0   &  0  &  0
  \end{array} \right], \]
where $K_1$ satisfy $J_0K_1 = - K_1J_0$ and $K_1(K_1 - J_0)J_0 = 0$; in other words,
$K_1$ are all the anti-commuting solutions of $J_0YJ_0 = YJ_0Y$.
\end{theorem}

\section{Numerical Examples}

\quad $\;$ To illustrate the main result of this paper, we present two concrete
Yang-Baxter-like matrix equations with a given singular constant matrix $A$.

\begin{example}
Let
\[A = \frac{1}{33}\left[ \begin{array}{rrrrrrrr}
		12  & 9   & -66 & 21  & 36  & 36  & 0   & 9   \\
		-12 & 57  & 0   & 12  & -36 & 30  & 33  & 57  \\
		64  & 4   & 88  & -64 & -72 & -72 & 0   & 4   \\
		2   & -4  & -55 & 31  & 39  & 6   & 0   & -4  \\
		90  & 18  & 99  & -90 & -93 & -93 & 0   & 18  \\
		-8  & 16  & 22  & 8   & -24 & 9   & 0   & 16  \\
		11  & -22 & -55 & -11 & 33  & 33  & -66 & 11  \\
		52  & -5  & 55  & -52 & -42 & -42 & -33 & -5     \end{array} \right].\]
Then, $A = WJW^{-1}$, where
\[	J = \left[ \begin{array}{rrrrrrrr}
		0 & 1 & 0 & 0 & 0 & 0 & 0 & 0 \\
		0 & 0 & 1 & 0 & 0 & 0 & 0 & 0 \\
		0 & 0 & 0 & 0 & 0 & 0 & 0 & 0 \\
		0 & 0 & 0 & 1 & 1 & 0 & 0 & 0 \\
		0 & 0 & 0 & 0 & 1 & 1 & 0 & 0 \\
		0 & 0 & 0 & 0 & 0 & 1 & 0 & 0 \\
		0 & 0 & 0 & 0 & 0 & 0 & -1 & 1 \\
		0 & 0 & 0 & 0 & 0 & 0 & 0 & -1
	\end{array} \right], \; \; W = \left[ \begin{array}{rrrrrrrr}
	1 & 2 & 2 & 1 & 1 & 1 & 0  & 0 \\
	1 & 1 & 0 & 0 & 2 & 2 & 1  & 0 \\
	4 & 0 & 1 & 0 & 0 & 0 & 0  & 0 \\
	0 & 0 & 0 & 1 & 0 & 0 & 0  & 0 \\
	6 & 0 & 1 & 0 & 1 & 0 & 0  & 0 \\
	0 & 0 & 2 & 0 & 0 & 1 & 0  & 0 \\
	0 & 0 & 1 & 0 & 0 & 0 & -1 & 1 \\
	3 & 0 & 0 & 0 & 0 & 0 & -1 & 0
	\end{array} \right].\]
	
By Theorem 3.9, all anti-commuting solutions of the simplified Yang-Baxter-like
matrix equation $JYJ = YJY$ are $K = \mbox{\rm diag}(K_{1}, 0, 0)$, where
 \[K_1 =\left[ \begin{array}{rrrrrrrr}
		0 &  y &  x  \\
		0 &  0 & -y  \\
		0 &  0 & 0
	\end{array} \right], \,  \forall \; x, y \in \mathbb{C}. \]
Hence, all anti-commuting solutions of the original Yang-Baxter-like
matrix equation $AXA = XAX$ are  $B = WKW^{-1} = (1/33) \times$
{\footnotesize
\[ \left[ \begin{array}{rrrrrrrr}
8y+4x   & 17y-8x & 44y-11x  & -8y-4x   & 12x-42y  & 12x-42y  & 0 & 17y-8x \\
12y+4x  & 9y-8x  & 33y-11x  & -12y-4x  & 12x-30y  & 12x-30y  & 0 & 9y-8x  \\
64y+16x & 4y-32x & 88y-44x  & -64y-16x & 48x-72y  & 48x-72y  & 0 & 4y-32x \\
0       &	0    & 0        & 0        & 0        & 0        & 0 & 0      \\
96y+24x & 6y-48x & 132y-66x & -96y-24x & 72x-108y & 72x-108y & 0 & 6y-48x \\
0       &	0    & 0        & 0        & 0        & 0        & 0 & 0      \\
0       &	0    & 0        & 0        & 0        & 0        & 0 & 0      \\
48y+12x & 3y-24x & 66y-33x  & -48y-12x & 36x-54y  & 36x-54y  & 0 & 3y-24x
\end{array} \right] \]}
\end{example}
\noindent with $x$ and $y$ being arbitrary complex numbers.

\begin{example}
Let \[J = \left[ \begin{array}{ccccccc}
		0 & 1 & 0 & 0 & 0 & 0 & 0  \\
		0 & 0 & 1 & 0 & 0 & 0 & 0  \\
		0 & 0 & 0 & 0 & 0 & 0 & 0  \\
		0 & 0 & 0 & 0 & 1 & 0 & 0  \\
		0 & 0 & 0 & 0 & 0 & 1 & 0  \\
		0 & 0 & 0 & 0 & 0 & 0 & 1  \\
		0 & 0 &0  & 0 & 0 & 0 & 0
	\end{array} \right] = \mbox{\rm diag}[J_3(0), J_4(0)]. 	\]
Partition a $7 \times 7$ matrix $K$ accordingly as
\[ K =  \left[ \begin{array}{cccc}
		K_1 & K_2 \\
		K_3 & K_4
	\end{array} \right]. \]
Then by Lemma 2.3, $JK = -KJ$ if and only if
\[  K_1 =	\left[ \begin{array}{rrr}
		k_{11} & k_{12}  & k_{13}  \\
		0      & -k_{11} & -k_{12} \\
		0      & 0       & k_{11}
	\end{array} \right], \; \;  K_2 =  \left[ \begin{array}{rrrr}
		0 & k_{22} & k_{23}  & k_{24}  \\
		0 & 0      & -k_{22} & -k_{23} \\
		0 & 0      & 0       & k_{22}
	\end{array} \right], \]
\[  K_3 =	\left[ \begin{array}{rrr}
		k_{31} & k_{32}  & k_{33}  \\
		0      & -k_{31} & -k_{32} \\
		0      & 0       & k_{31}  \\
		0      & 0       & 0
	\end{array} \right], \; \; K_4 =  \left[ \begin{array}{rrrr}
		k_{41} & k_{42}  & k_{43}  & k_{44}  \\
		0      & -k_{41} & -k_{42} & -k_{43} \\
		0      & 0       & k_{41}  & k_{42}  \\
		0      & 0       & 0       & -k_{41}
	\end{array} \right].\]
	
From Theorem 3.7 (ii), to satisfy the further equality $K(K - J)J = 0$,
the above sub-matrices $K_1, K_2, K_3, K_4$ of $K$ must satisfy

\begin{eqnarray*}
&  & \left[ \begin{array}{cc}
            K_1 & K_2 \\
            K_3 & K_4
            \end{array} \right] \left[ \begin{array}{cc}
            [K_1 - J_3(0)]J_3(0) & K_2J_4(0) \\
            K_3 J_3(0)           & [K_4 - J_4(0)]J_4(0)
            \end{array} \right] \\
& = & \left[ \begin{array}{cc}
K_1[K_1 - J_3(0)]J_3(0) + K_2K_3J_3(0)  & K_1K_2J_4(0) + K_2[K_4 - J_4(0)]J_4(0) \\
K_3[K_1 - J_3(0)]J_3(0) + K_4K_3J_3(0)  & K_3K_2J_4(0) + K_4[K_4 - J_4(0)]J_4(0)
            \end{array} \right] \\
& = & \left[ \begin{array}{cc}
            0 & 0 \\
            0 & 0
            \end{array} \right],
\end{eqnarray*}
from which
\begin{eqnarray}
&   \left\{    \begin{array}{ll}
			K_1 [K_1 - J_3(0)] J_3(0) + K_2 K_3J_3(0)  = 0,    \\
			K_1K_2J_4(0) + K_2[K_4 - J_4(0)]J_4(0) = 0, \\
			K_3[K_1 - J_3(0)]J_3(0) + K_4K_3J_3(0)  = 0,       \\
			K_3K_2J_4(0) + K_4[K_4 - J_4(0)]J_4(0) = 0.
		\end{array} \right. &
\end{eqnarray}
The above expressions of $K_1, K_2, K_3$, and $K_4$ lead to
\begin{eqnarray*}
&   & K_1 [K_1 - J_3(0)] J_3(0) + K_2 K_3 J_3(0) \\
& = & \left[ \begin{array}{rrr}
   	0 & k^2_{11} & -k_{11} -k_{22}k_{31} \\
   	0 & 0     & k^2_{11}                    \\
   	0 & 0       & 0
   \end{array} \right],
\end{eqnarray*}
\begin{eqnarray*}
&   & K_1 K_2 J_4(0) + K_2 [K_4 - J_4(0)] J_4(0) \\
& = & \left[ \begin{array}{rrrr}
	0 & 0 & -k_{22} k_{41} & -k_{12} k_{22} - k_{22} k_{42} -k_{22} + k_{23} k_{41} \\
	0 & 0 & 0         & -k_{22} k_{41}                                          \\
	0 & 0 &  0             & 0
             \end{array} \right],
\end{eqnarray*}
\begin{eqnarray*}
&   & K_3[K_1 - J_3(0)] J_3(0) + K_4 K_3 J_3(0) \\
& = & \left[ \begin{array}{rrr}
	0 & k_{41} k_{31} & k_{31} k_{12} - k_{31} + k_{41} k_{32} - k_{42} k_{31} \\
	0 & 0            & -k_{41} k_{31}                                        \\
	0 & 0            & 0                    \\
	0 & 0          & 0
	\end{array} \right],
\end{eqnarray*}
\begin{eqnarray*}
&   & K_3 K_2 J_4(0) + K_4[K_4 - J_4(0)] J_4(0) \\
& = & \left[ \begin{array}{rrrr}
  	0 & k^{2}_{41} & k_{31} k_{22} -k_{41} & k_{31} k_{23} - k_{32} k_{22}
        + 2k_{41} k_{43} - k^{2}_{42} - k_{42}                             \\
  	0 & 0          & k^{2}_{41}             & k_{31}k_{22} + k_{41}     \\
  	0 & 0          & 0                & k^{2}_{41}           \\
  	0 & 0          & 0             & 0
  \end{array} \right].
  \end{eqnarray*}




Therefore, (17) is equivalent to the system
\begin{eqnarray*}
	\left\{  \begin{array}{ll}
		k_{11} = 0,                                           \\
		k_{41} = 0,                                           \\
		k_{22} k_{31} = 0,                                    \\
		k_{22}(1 + k_{12} +	k_{42}) = 0,                      \\
		k_{31}(k_{12} - 1 - k_{42}) = 0,                      \\
		k_{31} k_{23}-k_{32} k_{22} -k_{42}^{2} - k_{42} = 0.
	         \end{array} \right.
\end{eqnarray*}
Solving the above, we obtain all anti-commuting solutions of (\ref{yb1}):
\[ K = \left[ \begin{array}{rrrrrrr}
		0 & k_{12} & k_{13}  & 0 & 0 & k_{23} & k_{24}  \\
		0 & 0      & -k_{12} & 0 & 0 & 0      & -k_{23} \\
		0 & 0      & 0       & 0 & 0 & 0      & 0       \\
		0 & k_{32} & k_{33}  & 0 & 0 & k_{43} & k_{44}  \\
		0 & 0      & -k_{32} & 0 & 0 & 0      & -k_{43} \\
		0 & 0      & 0       & 0 & 0 & 0      & 0       \\
		0 & 0     &0         & 0 & 0 & 0      & 0
	   \end{array} \right]; \]
\[ K = \left[ \begin{array}{rrrrrrr}
		0 & k_{12} & k_{13}  & 0 & 0  & k_{23} & k_{24}  \\
		0 & 0      & -k_{12} & 0 & 0  & 0      & -k_{23} \\
		0 & 0      & 0       & 0 & 0  & 0      & 0       \\
		0 & k_{32} & k_{33}  & 0 & -1 & k_{43} & k_{44}  \\
		0 & 0      & -k_{32} & 0 & 0  & 1      & -k_{43} \\
		0 & 0      & 0       & 0 & 0  & 0      & -1      \\
		0 & 0      & 0       & 0 & 0  & 0      & 0
	  \end{array} \right];  \]

\[ K = \left[ \begin{array}{rrrrrrr}
		0 & -(k_{42} + 1)                      & k_{13}                          & 0 & k_{22}
& k_{23}  & k_{24}  \\
		0 & 0                                  & k_{42} + 1                      & 0 & 0                & -k_{22} & -k_{23} \\
		0 & 0                                  & 0                               & 0 & 0                & 0       & k_{22}  \\
		0 & -\frac{k_{42}(k_{42} + 1)}{k_{22}} & k_{33}                          & 0 & k_{42}           & k_{43}  & k_{44}  \\
		0 & 0                                  & \frac{k_{42}(k_{42}+1)}{k_{22}} & 0 & 0                & -k_{42} & -k_{43} \\
		0 & 0                                  & 0                               & 0 & 0                & 0       & k_{42}  \\
		0 & 0                                  & 0                               & 0 & 0                & 0       & 0
	          \end{array} \right], \; \; k_{22} \neq 0; \]	

\[ K = \left[ \begin{array}{rrrrrrr}
0      & k_{42} + 1 & k_{13}        & 0 & 0      & \frac{k_{42}(k_{42} + 1)}{k_{31}} & k_{24}  \\
0      & 0          & -(k_{42} + 1) & 0 & 0      & 0                                 &
  -\frac{k_{42}(k_{42} + 1)}{k_{31}}  \\
0      & 0          & 0             & 0 & 0      & 0                                 & 0  \\
k_{31} & k_{32}     & k_{33}        & 0 & k_{42} & k_{43}                            & k_{44}  \\
0      & -k_{31}    & -k_{32}       & 0 & 0      & -k_{42}                           & -k_{43} \\
0      & 0          & k_{31}        & 0 & 0      & 0                                 & k_{42}  \\
0      & 0          & 0             & 0 & 0      & 0                                 & 0
	\end{array} \right], \; \; k_{31} \neq 0. \]
\end{example}

\section{Conclusions}

\ $\;$ We have found all anti-commuting solutions of the quadratic matrix
equation (\ref{yb}) with an arbitrary given matrix $A$. Our main result indicates
that, it is the Yang-Baxter-like matrix equation associated to the eigenvalue
$0$ of $A$ that leads to all such solutions.

The study of this problem was motivated in the continued effort to find all
solutions of (\ref{yb}), and our approach was based on a general study of some special
type Sylvester equation and an equivalent formation of (\ref{yb}) via Lemma 3.1. The
resulting solution process is divided into two steps. The first step was to solve
the anti-commutability equation $JY = -YJ$ and the second one focused on the equation
$Y(Y-J)J = 0$ in which $J$ is the Jordan canonical form of $A$. It turned out that
eventually, solving $JY = -YJ$ and $Y(Y - J)J = 0$ is essentially the same as solving
$J(0)Y = -YJ(0)$ and $Y[Y - J(0)]J(0)= 0$, as Theorem 3.9 has demonstrated.

We hope to generalize our result in the future to find other types of non-commuting
solutions of (\ref{yb}), and our goal is to find all solutions of the Yang-Baxter-like
matrix equation.
\section*{Funding}

\noindent This research is supported by the National Natural Science Foundation of China (12471133 and 11771378).

\section*{Disclosure statement}

No potential conflict of interest was reported by the authors.

\end{document}